\documentclass[12pt]{article}
\usepackage{amsxtra,amssymb,amsthm,amsmath,latexsym}

\theoremstyle{plain}
\numberwithin{equation}{section}

\def\R{{\mathbb R}}

\def\oH{{\overset{\circ}{H}}}
\def\oH1{{\overset{\circ}{H}\kern-.02in{}^1}}

\def\bee{\begin{equation*}}
\def\eee{\end{equation*}}
\def\be{\begin{equation}}
\def\ee{\end{equation}}

\begin{document}

Journ. Appl. Analysis, (2024),,,,, DOI: 10.1515/jaa-2024-0003

%\begin{titlepage}
\title{ On the Laplace transform
}

\author{Alexander G. Ramm\\
 Department  of Mathematics, Kansas State University, \\
 Manhattan, KS 66506, USA\\
ramm@ksu.edu\\
%,\\ %fax 785-532-0546, tel. 785-532-0580}
http://www.math.ksu.edu/\,$\sim$\,ramm}
%}

\date{}
\maketitle\thispagestyle{empty}

%%%%%%%%%%%%%%%%%%%%%%%%%%%%%%%%%%%%%%%%%%%%%%%%%
\begin{abstract}
\footnote{MSC: 44A10 }
\footnote{Key words:  Laplace transform; tempered distributions; analytic functions.
 }
Sufficient conditions are given for a function $F(p)$ to be the Laplace transform of a function $f(t)$ or a distribution $f$. No assumption on $f$ is given a priori.
It is not even assumed that $f=0$ for $t<0$.  

\end{abstract}

\section{Introduction}\label{S:1}

In reference book \cite{BE} one finds tables of the Laplace transform and \cite{S} is a text on this transform. In \cite{H} the theory of the Fourier transform can be found. 
It is often important to know if a function $F(p)$ is a Laplace transform of a function or a distribution, see, for example, \cite{R723}. In the literature it is usually assumed that $f(t)$
satifies some assumptions, for example, $f(t)=0$ for $t<0$, $|f(t)| < ce^{at}$, $a=const\ge 0$, see, for example, \cite{S}. By $c>0$ various constants are denoted. In this paper we do not make any assumptions on $f$.

Let us formulate our main result.

{\bf Theorem 1.}  {\bf Assume that $F(p), \, p=\sigma +is,$ is an analytic function in the half-plane $\sigma>\sigma_0$, there exist the limits $\lim_{\sigma\to \sigma_0} F(p)=F(\sigma_0+is)$ for almost every $s\in \R:=(-\infty, \infty)$ and 
\be\label{e1}
\lim_{|p_n|\to \infty} F(p_n)=0, 
\ee  
where $\lim_{n\to \infty}|p_n|= \infty$. 

Then
\be\label{e2}
F(p)=\int_0^\infty e^{-pt}f(t)dt,
\ee
and
\be\label{e3}
f(t)=\lim_{n\to \infty}\frac1{2\pi i}\int_{K_{\sigma_0 }}e^{qt}F(q)dq,
\ee
where $K_{\sigma n}:=(\sigma-in, \sigma+in),\, \sigma\ge \sigma_0\ge 0$.
 }

{\bf Proof.} By the analyticity of $F(p)$ and the Cauchy formula one gets
\be\label{e4}
F(p)=\frac 1 {2\pi i}\int_{L_{\sigma \,  n}}\frac {F(q)dq}{q-p},
\ee
where $p$ belongs to the domain with the closed boundary $L_{\sigma\, n}=-K_{\sigma_0\, n}\cup C_n$, $C_n=\sigma_0+ne^{i\phi}$, $-\frac \pi 2\le \phi\le \frac \pi 2$. The contour $L_{\sigma \,n}$ is oriented counterclockwise.

Let $n\to \infty$. Then 
\be\label{e5}
\lim_{n\to \infty}|\int_{C_n}\frac {F(q)dq}{q-p}|\le \lim_{n\to \infty}\frac n{n-|p|} \sup_{q\in C_n} |F(q)|=0,
\ee
because of our assumption \eqref{e1}. From \eqref{e5} and \eqref{e2} it follows that
\be\label{e6}
F(p)=\frac 1 {2\pi i}\int_{-K_{\sigma_0}}\frac {F(q)dq}{q-p}, \qquad -K_{\sigma_0}:=(\sigma_0+i\infty, \sigma_0- i\infty).
\ee
One has
\be\label{e7}
\frac 1 {p-q}=\int_0^\infty e^{-(p-q)t}dt, \qquad  Re p>Re q.
\ee
From \eqref{e6} and \eqref{e7} one gets
\be\label{e8}
F(p)=\int_0^\infty e^{-pt} \Big(\frac 1 {2\pi i}\int_{K_{\sigma_0}} F(q)e^{qt}dq\Big) dt=\int_0^\infty e^{-pt}f(t)dt,
\ee
where
\be\label{e9} 
f(t)=\frac 1 {2\pi i}\int_{K_{\sigma_{0}}} F(q)e^{qt}dq,
\ee
and
\be\label{e10} 
F(p)=\int_0^\infty e^{-pt} f(t)dt.
\ee
Let us prove that $f(t)$ in formula \eqref{e9} has the property $f(t)=0$ for $t<0$. 

Consider the identity which follows from the analyticity of $F(p)$ and of $e^{-q|t|}$ in the half-plane $\sigma>0$ and the decay of $e^{-q|t|}$ when $Re\, q\to \infty$ and $|t|>0$: 
\be\label{e11} 
\int_{L_{\sigma_0 n}}e^{-q|t|}F(q)dq=0.
\ee
We have
\be\label{e12} 
\lim_{n\to \infty} |\int_{C_n}e^{-q|t|}F(q)dq|\le \lim_{n\to \infty}\sup_{q\in C_n}|F(q)| \int_{-\frac \pi 2}^{\frac \pi 2} e^{-n\cos \phi |t|}n d\phi=0,
\ee
because 
\be\label{e13}  
 \int_{-\frac \pi 2}^{\frac \pi 2} e^{-n\cos \phi |t|}n d\phi=2\int_0^{\frac \pi 2} e^{-n \cos \phi |t|}n d\phi=2\int_0^{\frac \pi 2}e^{-n|t| \sin\phi} d\phi=
O(\frac 1 {n|t|}), \quad n\to \infty.
\ee
This estimate  follows from the known inequality $\sin \phi\ge \frac 2 \pi \phi$ and a simple estimate:
\be\label{e14}  
 \int_0^{\frac \pi 2}e^{-\frac 2 \phi n|t| \sin\phi} d\phi=
O(\frac 1 {n|t|}), \quad n\to \infty.
\ee
From \eqref{e11}, \eqref{e12} and  \eqref{e14}  it follows that
\be\label{e15} 
\frac 1 {2\pi i}\int_{K_{\sigma_0}}e^{-q|t|}F(q)dq=0, \quad |t|>0.
\ee
This means that $f(t)=0$ for $t<0$.

Theorem 1 is proved.  \hfill$\Box$

 \section{Various comments}\label{S:2}

  a)   Conditions for $F(p)$ to be a Laplace transform are practically important, see, for example, \cite{R723}, where the equation $h(t)=g(t)+\int_0^\infty (t-s)^{-\frac 5 4}h(s)ds$ is studied. The integral is hypersingular, it diverges classically. In   \cite{R723}, pp.23-31, this classically divergent integral is defined,
the integral equation is solved explicitly, the Laplace solution of this equation is $L(h)=\frac{L(f)}{1+c_1p^{\frac 1 4}}$, where $c_1:=|\Gamma(-\frac 1 4)|>0$.
If $f=f(t)$ is smooth and rapidly decaying when $t\to \infty$, then $|L(h)|\le \frac c {1+|p|}, \, |p|\gg 1,\,\, Re\, p\ge 0$. One proves that this $L(h)$
satisfies assumptions of Theorem 1, $L(h)$ is analytic in the half-plane $\sigma>0$, so it is indeed a Laplace transform of a function. One also sees that 
$\lim_{|p|\to\infty} p^{\frac 54}L(h)=0$. Therefore, $f(0)=0$, see Lemma 1 below.

The question of practical interest is: can one prove that the solution $h=h(t)$ is a uniformly bounded function on $\R_+$ and $h(0)=0$? Theorem 1 and 
Lemma 1, proved below in item c), allow to answer both questions.

b)  Assume that $\sigma_0=0$ and $F(is)\in L^\ell(\R):=L_\ell$, $\ell\ge 1$, where $L_\ell$ are the Lebesgue spaces,  
\be\label{e16}
\frac 1 {2\pi}\int_{-\infty}^\infty  |F(is)|^\ell ds<\infty.
\ee
If $\ell\in [1, 2]$, then the Fourier transform $f(t)=\frac 1 {2\pi}\int_{-\infty}^\infty  F(is)e^{-ist}ds$ belongs to $L_{\frac \ell{\ell-1}}$, see \cite{H}, p.165.
One has the  Hausdorff-Young estimate \cite{H}, p.165:
\be\label{e17}
\|f(t)\|_{L_{\frac \ell{\ell-1}}}\le (2\pi)^{\frac {\ell-1}{\ell}}\|F\|_{L_\ell}.
\ee
Therefore, if $\ell\in [1,2]$, then $f(t)$ is a function. If $\ell>2$, then $f$ may be a tempered distribution, see \cite{H}, pp.163-164.

c) {\bf Example.}  Let $F(p)=\frac 1{(1+p)p^{\frac 1 4}}, \,p=\sigma+is, \, \sigma>0 $. This function is analytic in the half-plane $\sigma>0$ and all the
conditions of Theorem 1 are satisfied. Therefore, this function is a Laplace transform of $f$, $F(p)=L(f)$. Since $F(is)\in L_1$, the function $f(t)$ is uniformly bounded.

Since $\lim_{|p|\to \infty} |p|^{b} F(p)=0$, $b=\frac 5 4>1$, we conclude that $f(0)=0$.

This result is a consequence of the following Lemma 1 in which the conditions of Theorem 1 are assumed to hold without explicitly repeating them and it is  also assumed
that $\sigma_0=0$ and $F(is)\in L_1$.

{\bf Lemma 1. If $\sigma>0$, $b>1$, $F(p)$ satisfies assumptions of Theorem 1 and $\lim_{|p|\to \infty} |p|^b |F(p)|=0$, then $f(0)=0$. }

{\bf First proof of Lemma 1.} By Theorem 1, we have formula \eqref{e10}.  Write this formula as
\be\label{e18}
F(p)=(\int_0^\epsilon +\int_\epsilon^\infty) e^{-pt}f(t)dt:=I_1+I_2,
\ee
where $\epsilon>0$ is a small number which will be chosen later.  Let $p>0$ tend to $\infty$. Since $f(t)$ is uniformly bounded, one has:
\be\label{e19} 
|I_2|\le c\frac{e^{-\epsilon  p}}p,
\ee
and
\be\label{e20} 
|I_1|\le c\frac{1-e^{-\epsilon  p}}p f(0)(1+O(\epsilon p)),
\ee
because for small $\epsilon$, continuous $f(t)$ and $f(0)\neq 0$ one has $f(t)=f(0)(1+O(\epsilon p))$.

Choose $\epsilon=\frac 1 {p^{1/2}}, \,p\to +\infty$. Then $\lim_{p\to \infty}|p|^bI_2=0$. Since  Since $b>1$,  we have $\lim_{p\to \infty}|p|^bF(p)=0$.
Therefore,  $\lim_{p\to \infty}|p|^bI_1=0$. This can happen only if $f(0)=0$, see \eqref{e20}.

Lemma 1 is proved.

{\bf Second proof of Lemma 1.} By formula \eqref{e9} with $\sigma_0=0$ and $t=0$ one has
\be\label{e21} 
f(0)=\frac 1{2\pi}\int_{-\infty}^\infty F(q)ds=0.
\ee
Indeed, 
\be\label{e22} 
\lim_{n\to \infty}\int_{C_n}F(q)dq=0,
\ee
 because $b>1$. By analyticity of $F$ one has
\be\label{e23} 
\int_{L_n}F(q)dq=0.
\ee
From \eqref{e22} and \eqref{e23} passing to the limit $n\to \infty$ one derives formula \eqref{e21}.

The second proof of Lemma 1 is complete. \hfill$\Box$

d)  Let $\R_+=(0,\infty),\,\R=(-\infty, \infty)$, $p=\sigma +is$.  Consider the following result of the type of a Paley-Wiener (PW) theorem, see \cite{R}, p.405.

{\bf Proposition 1.  Suppose $F=F(p)$ is given by formula \eqref{e2} with $f\in L^2(\R_+)$. Then $F$ is analytic in the half-plane $\sigma>0$ and
$\sup_{\sigma\ge 0}\int_{-\infty}^{\infty}|F(p)|^2ds\le c $. 

Conversely, if $F$ is analytic in the half-plane $\sigma>0$ and
$\sup_{\sigma\ge 0}\int_{-\infty}^{\infty}|F(p)|^2ds\le c $, then $F$ is  given by formula \eqref{e2} with $f\in L^2(\R_+)$.}

{\bf Proof of Proposition 1.} Let $p=\sigma+is=r(\cos \phi +i\sin \phi), \, r=|p|$. For any fixed $\phi\in (-\frac \pi 2, \frac \pi 2)$ one has by the Cauchy inequality:
\be\label{e24}
|\int_0^\infty e^{-pt}f(t)dt|^2\le \int_0^\infty |f(t)|^2dt \int_0^\infty e^{-2r\cos \phi\,  t}=  \int_0^\infty |f(t)|^2dt \frac 1 {2r\cos 
\phi}.
\ee
Thus, $\lim_{|p|\to \infty}F(p)=0$ for any $\phi$ not equal to $\pm \frac \pi 2$. For $\phi =\pm \frac \pi 2$,  one has $F(p)=\int_0^\infty e^{-ist}f(t)dt:=\mathcal{F}(f)$, where  $\mathcal{F}$ is the Fourier transform. Since $f\in L^2(\R_+)$, there is a sequence $s_j\to \infty$ such that
$\lim_{j\to \infty}F(is_j)=0$. By the Parseval identity one has:
\be\label{e25}
\sup_{\sigma\ge 0}\int_{-\infty}^{\infty}|F(\sigma +is)|^2ds=2\pi \sup_{\sigma\ge 0}\int_0^\infty |e^{-\sigma t}f(t)|^2dt\le 2\pi  \int_0^\infty |f(t)|^2dt\le c.
\ee
Therefore, the first part of Proposition 1 is proved. 

To prove the second part, assume that $F$ is analytic  in the half-plane $\sigma>0$ and 
\be\label{e26}
\sup_{\sigma\ge 0}\int_{-\infty}^{\infty}|F(p)|^2ds\le c.
\ee
From  formula \eqref{e26} it follows that the limit $\lim_{\sigma\to 0}F(p)=F(is)$ exists a.e. (for almost every $s\in \R$) and condition \eqref{e1} holds. By Theorem 1 formula \eqref{e2} holds and by formula \eqref{e25}  $ \int_0^\infty |f(t)|^2dt\le c.$

Proposition 1 is proved. \hfill$\Box$

e)  Assume that $F(p)=e^{\frac 1 {p^2}}$. Is $F(p)$ a Laplace transform of a function? Let us prove that the function  $e^{\frac 1 {p^2}}$ is not a Laplace transform of any function $f(t)$ such that $f=0$ for $t<0$ and $\int_0^\ell |f(t)|^\ell dt<c$, $\ell\ge 1$.  Let $p=\sigma+is=\sigma>0$.  
Suppose to the contrary that $e^{\frac 1 {\sigma^2}}=\int_0^\infty e^{-\sigma t}f(t)dt$.  We have
\be\label{e27}
e^{\frac 1 {\sigma^2}}=\int_0^\infty e^{-\sigma t}f(t)dt\le (\int_0^\infty|f(t)|^\ell dt)^{\frac 1 \ell}(\int_0^\infty e^{-\ell' \sigma t}dt)^{\frac 1 {\ell'}}\le 
 \frac c{(\ell' \sigma)^{\frac 1 {\ell'}}},
\ee
where the H\"{o}lder inequality is used and $\ell'$ is defined by the formula $\frac 1 \ell + \frac 1{ \ell'}=1$.
 When $\sigma\to 0$ the right side of \eqref{e27} tends to infinity less fast than the left side.  Therefore, the first equality sign in equation \eqref{e27} is impossible.

  \section{Conclusion}\label{S:3}
Sufficient conditions are given for a function $F(p)$ to be the Laplace transform of $f(t)$. No assumptions about $f(t)$ are imposed.

\section{Conflict of interest}\label{S:4}
There is no conflict of interest. 
%\newpage 

\end{document}